\documentclass[journal]{IEEEtran}
\usepackage{float}
\usepackage{makeidx}  % allows index generation
\usepackage{graphicx} % standard LaTeX graphics tool
\usepackage{mathrsfs}
\usepackage{arydshln}
\usepackage[tbtags]{amsmath}         % added by editor
\usepackage{amssymb}                 % added by editor
\usepackage{epsfig}                  % added by editor
\usepackage{latexsym}                % added by editor
\usepackage{array}                  % added by editor
\usepackage{ifthen}
\usepackage{epstopdf}
\usepackage{booktabs}
\usepackage{subfigure}
\usepackage{theorem}
\usepackage{color}
\usepackage{comment}
\usepackage{CJK}
\usepackage{epstopdf}
\usepackage{cases}
\usepackage{subfigure}

\newtheorem{theorem}{Theorem}
\newtheorem{proof}{Proof}
\newtheorem{lemma}{Lemma}

\newtheorem{assum}{Assumption}

\newtheorem{remark}{Remark}

\newcommand{\be}{\begin{equation}}
\newcommand{\ee}{\end{equation}}

\ifCLASSINFOpdf
\else
\fi
\begin{document}
%
% paper title
% Titles are generally capitalized except for words such as a, an, and, as,
% at, but, by, for, in, nor, of, on, or, the, to and up, which are usually
% not capitalized unless they are the first or last word of the title.
% Linebreaks \\ can be used within to get better formatting as desired.
% Do not put math or special symbols in the title.
\title{Platooning of Connected Vehicles with Directed
Graph: $H_\infty$ Robustness Analysis and Synthesis}
%
%
% author names and IEEE memberships
% note positions of commas and nonbreaking spaces ( ~ ) LaTeX will not break
% a structure at a ~ so this keeps an author's name from being broken across
% two lines.
% use \thanks{} to gain access to the first footnote area
% a separate \thanks must be used for each paragraph as LaTeX2e's \thanks
% was not built to handle multiple paragraphs
%

\author{Jilie~Zhang, Tao~Feng*%, Yin Tong~\IEEEmembership{Member,~IEEE}, Zhanbo Sun and Kehu Luo
        %and~Xiaomin~Wang <-this % stops a space
\thanks{Jilie~Zhang and Tao Feng are with the School of Information Science and Technology,
Southwest Jiaotong University, Chengdu, Sichuan, 610031, P.R. China. (e-mail: jilie0226@163.com; sunnyfengtao@163.com; yintong@swjtu.edu.cn)}% <-this % stops a space

%\thanks{Huaguang Zhang is with the School of Information Science and Engineering, Northeastern University, Shenyang, 110819, China. (e-mail: hgzhang@ieee.org)}
\thanks{This work was supported by the National Natural Science Foundation of China under grant 61803317 and 61603314.}% <-this % stops a space
%\thanks{Manuscript received April 19, 2020; revised August 26, 2020.}
}

% note the % following the last \IEEEmembership and also \thanks -
% these prevent an unwanted space from occurring between the last author name
% and the end of the author line. i.e., if you had this:
%
% \author{....lastname \thanks{...} \thanks{...} }
%                     ^------------^------------^----Do not want these spaces!
%
% a space would be appended to the last name and could cause every name on that
% line to be shifted left slightly. This is one of those "LaTeX things". For
% instance, "\textbf{A} \textbf{B}" will typeset as "A B" not "AB". To get
% "AB" then you have to do: "\textbf{A}\textbf{B}"
% \thanks is no different in this regard, so shield the last } of each \thanks
% that ends a line with a % and do not let a space in before the next \thanks.
% Spaces after \IEEEmembership other than the last one are OK (and needed) as
% you are supposed to have spaces between the names. For what it is worth,
% this is a minor point as most people would not even notice if the said evil
% space somehow managed to creep in.

% The paper headers
\markboth{%Journal of \LaTeX\ Class Files,~Vol.~XX, No.~X, August~XXXX
}%
{Shell \MakeLowercase{\textit{et al.}}: Bare Demo of IEEEtran.cls for IEEE Journals}
% The only time the second header will appear is for the odd numbered pages
% after the title page when using the twoside option.
%
% *** Note that you probably will NOT want to include the author's ***
% *** name in the headers of peer review papers.                   ***
% You can use \ifCLASSOPTIONpeerreview for conditional compilation here if
% you desire.

% If you want to put a publisher's ID mark on the page you can do it like
% this:
%\IEEEpubid{0000--0000/00\$00.00~\copyright~2015 IEEE}
% Remember, if you use this you must call \IEEEpubidadjcol in the second
% column for its text to clear the IEEEpubid mark.

% use for special paper notices
%\IEEEspecialpapernotice{(Invited Paper)}

% make the title area
\maketitle

% As a general rule, do not put math, special symbols or citations
% in the abstract or keywords.
\begin{abstract}
This paper revisits the distributed ${ H}_\infty$ robustness controller design for the platooning of connected vehicles.
Recently, the relevant result subjected to the undirected topology has been studied,
in the light of the symmetry of Laplace matrix.
It is well known that the same problem is more challenging for the \emph{directed} topology,
since the Laplace matrix ceases to be symmetric. In this paper, the problem is solved by introducing more weighting parameters and setting suitable values for them.
Then we show that the introduced weighting parameters lead to a positive effect on robustness, and solve the problem of feedback high gain.
Finally, two numerical simulations and a practical simulation based on Next Generation Simulation (NGSIM) dataset are applied to illustrate the effectiveness of our method.
\end{abstract}

% Note that keywords are not normally used for peerreview papers.

% For peer review papers, you can put extra information on the cover
% page as needed:
% \ifCLASSOPTIONpeerreview
% \begin{center} \bfseries EDICS Category: 3-BBND \end{center}
% \fi
%
% For peerreview papers, this IEEEtran command inserts a page break and
% creates the second title. It will be ignored for other modes.
\IEEEpeerreviewmaketitle

\section{Introduction}

\IEEEPARstart{T}{he} platooning of connected vehicles has been widely investigated for the past decades.
It gets increasing attention in the control community for its applications, such as automated highway systems (AHS),
intelligent transportation system (ITS), etc.
For a comprehensive literature review, readers may refer to some recent survey papers \cite{OPA2015,bpcest2012, t2020, wbswlb2020, wtm2009} and references therein.

In the platooning of connected vehicles,
a fundamental problem is inter-vehicle distance control.
In this scenario, the platoon control can ensure that vehicles are driving with the pre-specified space among them,
at the desired speeds \cite{gllg2019, kpd2017}.
{However, in the real world, the dynamics of a platoon system
%\textcolor[rgb]{may well be1.00,0.00,0.00}{influenced by the aging of vehicles (e.g., aging of electric component or abrasion of mechanical structure, etc.).}
%\textcolor[rgb]{1.00,0.00,0.00}{In addition, the platoon system is also}
is inevitably affected by the environment (e.g., aerodynamic drag,  frictional drag and slope, etc.).
The cases will change the dynamic of the platoon system, even may well result in that the system be unstable under the original controller.}

In a word, the robustness problem of platoon control systems becomes particularly important.
{However such a problem has been solved by some existing works \cite{Ploeg2014,lglwyc2018, zllr2018}.}
The paper \cite{Ploeg2014} synthesises an $H_{\infty}$ controller to achieve the string stability for the platooning with linear dynamics, by solving a linear
inequality matrix.
The paper \cite{lglwyc2018} designs the distributed ${H}_\infty$ control and analyzes the robustness in the frequency domain,
addressing the dynamics with an uncertainty model in a platoon system.
While the paper \cite{zllr2018} deals with the scenario where the platoon system is with the external disturbance in the time domain.
However both of them address the case that the interaction topologies of information are undirected among vehicles.
This case is relatively easy to be solved, since the undirected topology's Laplace matrix is symmetric. In other words, it is diagonalizable. If the matrix involving to graph is diagonalizable is decisive to analyze the robustness problem.
The platooning of vehicles over the undirected topology is more robust than the directed case in general, but it is at the expense of more sensors and network bandwidth in practice.
Intuitively, it implies that a system is more robust if it can get more useful information.
However, the platoon problem over directed graph is more important and more consistent with reality, taking into account the cost of devices and network bandwidth. At present papers \cite{zllbh2017} and \cite{zlwcl2016} have been devoted to solving it.
However all works subjected to directed topology are usually coupling, thus the corresponding robustness analysis is always a great challenge.

Motivated by the aforementioned work \cite{zllr2018},
this paper aims to propose a novel distributed ${ H}_\infty$ controller by introducing more weighting parameters,
which not only solves the inter-vehicle distance control problem over the directed topology, with the external disturbance, but also makes the robustness analysis in the time domain possible.
More specifically, we intend to change the weights of the graph in disguise and make the collected system decoupled by designing suitable parameters.
Compared with the work in \cite{zllr2018}, the contributions of our work are listed as follows:
\begin{enumerate}
  \item The robustness analysis is solved for the directed topology scenario.
        The difficulty lies in that the Laplace matrix ceases to be diagonalizable (it is hard to be decoupled).
      As far as we know, no related results yet emerge at present;
  \item It is shown that the robustness of platoon control is also dependent on the introduced weighting parameters, except for the number of the vehicles $N$ and the topological structure;
  \item
{The feedback high gain problem is circumvented by our method.
   In general, the control signal is bounded and it may not be infinitely large. If the coupling strength is too large, it is difficult to implement in practice.} %In addition, the fuel consumption and stationarity can be also improved in a platooning.
\end{enumerate}

%Note that, to make a sharp contrast with \cite{zllr2018}, here we use the similar expression to illustrate our approach as much as possible,
%even structure and title of the article.

The rest of this paper is organized as follows. Section \ref{s2} formulates the problem to be solved in this paper.
The robustness analysis with the directed topology is provided in Section \ref{s3}.
The distributed ${ H}_\infty$ synthesis for the platooning is shown in Section \ref{s4}.
Numerical examples in Section \ref{s5} illustrate the effectiveness of the proposed controller. Section \ref{s6} concludes the paper.
%\hfill mds

%\hfillAugust 26, 2015

\section{Problem Statement}\label{s2}
Here we consider the case as the same as that in \cite{zllr2018}. The platoon problem of connected $N + 1$
homogeneous vehicles is studied, with a leading node indexed by $0$.
The remaining vehicles are indexed from $1$ to $N$, which are referred to following nodes.
Our objective is to analyze the robustness and synthesize a distributed $H_\infty$ controller to ensure that all vehicles move at the desired speed,
meanwhile maintain the specific distances.
The only difference is that the information flow among vehicles is directed.
It is more challenging for this case, which is not yet solved up to now.

Here, we use the following dynamic model to express the acceleration response
of vehicle dynamics
\begin{align}
\tau \dot{a}_i (t) + a_i (t) = u_i (t) + w_i (t),
\end{align}
where $a_i (t)$ denotes the acceleration of vehicle $i$; $\tau$ is the
time delay parameter; $u_i$ is the control input, and $w_i (t)$ is the
disturbance from the exosystem, such as the aerodynamic drag, frictional drag and slope, etc.

Letting $p_i (t)$, $v_i (t)$ and $a_i (t)$ to denote respectively each vehicle's position, velocity and
acceleration, a state-space system of a vehicle dynamics is written as
\begin{align} \label{e2}
\dot{x}_i (t) = Ax_i (t) + B_1u_i (t) + B_2w_i (t),
\end{align}
where
\[
x_i =\left[
      \begin{array}{c}
        p_i(t) \\
        v_i(t)  \\
        a_i(t)  \\
      \end{array}
    \right],
A = \left[
      \begin{array}{ccc}
        0 & 1 & 0 \\
        0 & 0 & 1 \\
        0 & 0 & -\frac{1}{\tau} \\
      \end{array}
    \right],
B_1 = B_2 = \left[
      \begin{array}{c}
        0 \\
        0 \\
        \frac{1}{\tau}\\
      \end{array}
    \right].
\]
{The dynamic system (\ref{e2}) is a simplified model of vehicle with parasitic time delays and lags (see \cite{Darbha1994, Xiao2011} and \cite{SWAROOP1994}, etc.)}.

{In a platooning, the allowable communication connections subjected to directed graphs among vehicles is a more general case.} %For more comprehensive descriptions on graph theory, please see \cite{gr2013} and the references therein.

First of all, a neighbor set of node $i$ is defined as
\[\mathbb{N}_i=\{j\in {\cal V}_{N}|a_{ij}=1\},\]%is the Laplace matrix of the directed graph..
where ${\cal V}_{N}$ is a set of $N$ nodes, i.e., ${\cal V}_{N} = \{1, 2,\ldots , N\}$;
$a_{i j} = 1$ is the weight from neighbour node $j$ to node $i$ in a directed graph.

To model the communications from the leader to the
followers, we define an augmented directed graph $\mathcal{G}_{N+1} =
(\mathcal{V}_{N+1}, \mathcal{E}_{N+1})$ with a set of $N + 1$ vertices ($\mathcal{V}_{N+1} =
\{0, 1, \ldots, N\}$), which includes both the leader and the followers
in the platooning; $\mathcal{E}_{N+1}$ is a set of edges, which belongs to $\mathcal{V}_{N+1}\times \mathcal{V}_{N+1}$.
We use a pinning matrix $\mathcal{P} = \mbox{diag}\{g_1, g_2, . . . , g_N\}$ to denote how each follower connects
to the leader: $g_i = 1$ (or a positive constant) if $(0, i ) \in \mathcal{E}_{N+1}$, otherwise $g_i = 0$.
Note that the interactive information is directed among nodes.

\begin{assum}
The graph of a platooning $\mathcal{G}_{N+1}$ is a star topology \cite{zllr2018}.
%The graph of a platooning $\mathcal{G}_{N+1}$ contains at least
%one spanning tree \cite{Feng2021} rooting at the leader.
\end{assum}
The assumption means %every follower is globally reachable from the leader in $\mathcal{G}_{N+1}$.Namely
{the leader can send its information to every follower directly.}
Note that the commonly used topologies in the platooning belongs to the star topology, ( see Fig.1 in \cite{zlwcl2016}).

%, which is a prerequisite to
%guarantee the internal stability of a platoon.

\textbf{Objective:} The platoon control makes vehicles achieve the same
speed and the desired spacing between the front and the rear vehicles:
\begin{align}
&\lim_{t\rightarrow\infty}\|v_i(t)  - v_0(t) \|= 0; \nonumber\\
&\lim_{t\rightarrow\infty}\|p_i(t)  - p_{i-1}(t)  - d_{i,i-1}\|= 0;\ i = 1, \ldots , N,
\end{align}
where $d_{i,i-1}$ is the desired spacing between vehicles $i$ and
 $i-1$. $v_0$ is the leader's velocity.

In this paper, to analyze the robustness of a platooning over the directed graph,
we modify the local controller (8) in \cite{zllr2018} as
\begin{align}\label{e4}
u_i = -ck^T( g_i \hat{x}_i+\sum_{j\in \mathbb{N}_i}(d_i\hat{x}_i-d_{ij}\hat{x}_j) ), i = 1,\ldots, N,
\end{align}
where $k = [ k_p\quad k_v\quad k_a ]^T$, $c$ is the coupling strength,
and $g_i, d_i, d_{ij}>0$ are the weighting parameters, which are introduced deliberately to strengthen the self signal (by $g_i, d_i$),
and weaken the neighborhood information (by $d_{ij}$).
The aim is to implement the decoupling controller. {$\hat{x}_i = [\hat{p}_i \quad \hat{v}_i \quad  \ \hat{a}_i]^T$ is the tracking error for vehicle $i$, that is,
\begin{align*}
\hat{p}_i&=p_i(t) -p_0(t) -d_{i,0}, \\
\hat{v}_i&=v_i(t) -v_0(t) ,\\
\hat{a}_i&=a_i(t) -a_0(t) .
\end{align*}
where, $\hat{p}_i$  is the spacing error between the vehicles $i$ and $0$ with $d_{i,0}$ which is the desired spacing between the vehicles $i$  and $0$.
$\hat{v}_i$ and $\hat{a}_i$ are respectively the velocity error and acceleration error between the vehicles $i$ and $0$. %$\hat{a}_i$  is the acceleration error between the vehicle $i$  and the vehicle $0$.}
\begin{remark}
Note that in the controller (\ref{e4}), we use the absolute error of neighbours' information by installing some amplifiers with the gain $d_i$ and $d_{ij}$, rather than a relative error.
However the controller with relative error is  more widely used in the platoon control community. For this purpose, we can rewrite the controller (\ref{e4}) as
\begin{align}\label{e444}
u_i = -ck^T(  g_i \hat{x}_i+\sum_{j\in \mathbb{N}_i} \tilde{d}_i\hat{x}_i+\sum_{j\in \mathbb{N}_i}d_{ij}\hat{e}_{ij} ), i = 1,\ldots, N,
\end{align}
where $\hat{e}_{ij}=\hat{x}_i-\hat{x}_j$, and $\tilde{d}_i=d_i-d_{ij}$. Obviously, a vehicle only needs its own feedback, as the relative error can be obtained.
In fact, the controller (\ref{e4}) with absolute error can be unified as the relative one (\ref{e444}). $\hfill\blacksquare$
\end{remark}
The collective formulation of the controller (\ref{e4}) is
\begin{align}\label{l4}
U = -c\mathcal{M}\otimes k^TX,
\end{align}
where $X=[\hat{x}_1^T\quad \ldots\quad\hat{x}_N^T]^T$, $U=[u_1\quad \ldots\quad u_N]^T$, $\mathcal{M}={\cal L}_d+ \mathcal{P}$,
${\cal L}_d= {\cal D}-{\cal A}$, ${\cal D}=\mbox{diag}\{d_i \sum^N_{j=1}a_{ij}\}$ and ${\cal A}=[d_{ij}a_{ij}]_{N\times N}$.

Then, we can write the closed-loop dynamics of the platooning as
\begin{align}\label{e6}
\dot{X} &= A_c X+ BW,\nonumber\\
Y &= CX,
\end{align}
where $W=[w_1\quad \ldots\quad w_N]^T$, $A_c = I_N \otimes A - c\mathcal{M}\otimes B_1k^T $, $B=I_N\otimes B_2$, $C = I_N \otimes C_1$, $C_1 = [1, 0, 0]$, and
$Y$ is the tracking error of positions (i.e., $Y= [ \hat{p}_1\quad \hat{p}_2\quad \ldots\quad \hat{p}_N ]^T$), which denotes the output of a
platooning. Under the zero initial
tracking errors, the transfer function from the disturbance $W$ to the position $Y$ is
\begin{align}\label{e7}
G(s) &= C(s I_{3N} - A_c)^{-1}B,\nonumber\\
& = [I_N(\tau s^3+s^2) +c ({\cal L}_d + {\cal P})(k_p+k_vs+k_as^2)]^{-1}.
\end{align}
{Generally, for a directed graph it ceases to be trivial to analytically get the $H_\infty$ norms. The inverse operation in (\ref{e7}) is the greatest challenge, especially
simultaneously involving the factor ${\cal L}_d + {\cal P}$
(It ceases to be symmetric as ${\cal L} + {\cal P}$ in \cite{zllr2018}).
In the next section, we will focus on the analysis and synthesis of the distributed robustness controller
with a guaranteed performance by the novel controller (\ref{e4}).}

\section{Main Results on Robustness Analysis of Platooning with Directed Graph}\label{s3}
In the light of the idea in \cite{zllr2018}, if only ${\cal M}$ is diagonalizable, the robustness analysis and controller synthesis are feasible.
In the following lemma, we give a sufficient condition that the matrix ${\cal M}$ can be diagonalized.

{\begin{lemma}\label{l1}
For ${\cal M}$,
if a sequence $\{ \kappa_1, \cdots, \kappa_N\} =\cal N$  satisfies
\begin{align}
o_{\kappa_1}>r_{\kappa_1},\;
o_{\kappa_{i+1}}-o_{\kappa_{i}}>r_{\kappa_{i+1}}+r_{\kappa_{i}}, \; i=1,\ldots, N-1, \label{ddiag}
\end{align}
where $o_i=g_i+ d_i\sum_{j=1}^N a_{ij}$ and $r_i= \sum_{j=1}^N   d_{ij}a_{ij}$,  then ${\cal M}$ is diagonalizable.
\end{lemma}}

\begin{proof} We know that the diagonal entry of the $i$-th row of
${\cal M}$ is $o_i$ and the sum of the absolute values of the non-diagonal entries in the $i$-th row is
$r_i=\sum_{j=1}^N d_{ij}  a_{ij}$. It can define a Ger\v{s}gorin disc centred at $o_i$
with radius $r_i$, denoted as $D(o_i, r_i)$.

Under the condition \eqref{ddiag}, the Ger\v{s}gorin discs are
$D(o_{\kappa_i}, r_{\kappa_i})$ in the order of $i=1,\cdots, N$,
from the left to the right, and they do not intersect; see Fig. \ref{f22}.
Therefore, by Ger\v{s}gorin theorem, all eigenvalues of $ {\cal M}$ are distinct positive real numbers.
In other words,  ${\cal M}$ is diagonalizable. $\hfill\blacksquare$
\end{proof}

\begin{figure}[h!]
  \centering
  % Requires \usepackage{graphicx}
  \includegraphics[width=3in]{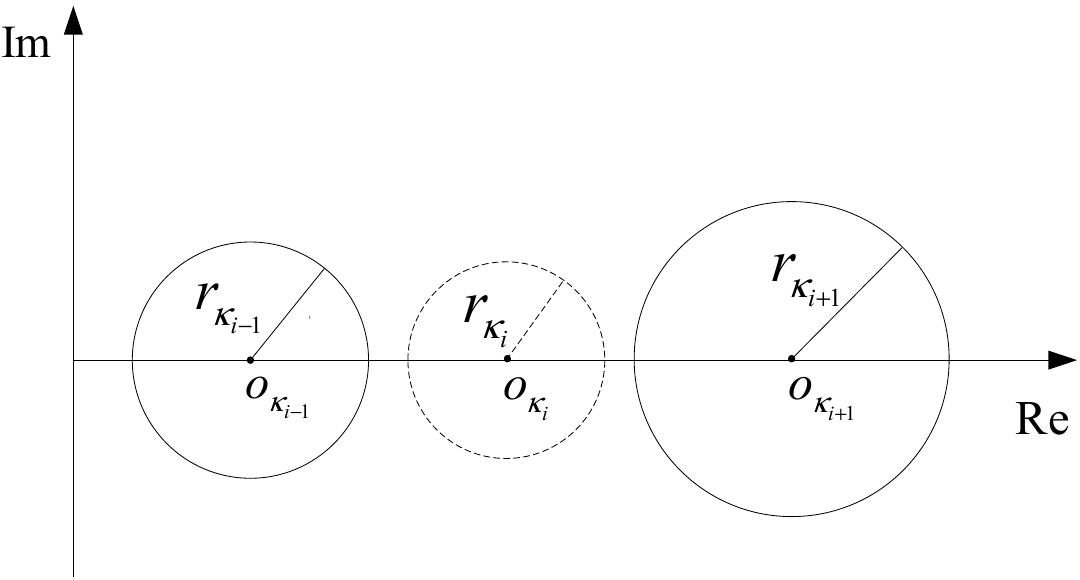}
  \caption{Non-intersect Ger\v{s}gorin discs.}\label{f22}
\end{figure}

%\vspace{1in}
\begin{remark}\label{r000}
The proof of Lemma 1 shows that if the center of the leftmost Ger\v{s}gorin
disc is sufficiently large with constant radius, then the diagonal entries of ${\cal M}$
are sufficiently large relative to the non-diagonal ones.
{In brief, we use the idea of diagonal dominance to make the diagonal elements sufficiently large,
or the non-diagonal elements sufficiently small, then
the matrix ${\cal M}$ approximates to a diagonal matrix.}
This property will be used in the proof of the following theorem.  $\hfill\blacksquare$
\end{remark}

{\begin{theorem}\label{t1}
 Consider a platooning, which consists of homogeneous vehicles, over the \textbf{directed
graph} given by (\ref{e7}). Under the conditions in Lemma \ref{l1}, using any stabilizing feedback
gains, the robustness measure $\gamma$-gain satisfies
\[
\gamma \leq \inf_{d_i,d_{ij}}\left\{\frac{1}{\lambda_{\min}k_p}\sqrt{\frac{\lambda_{max}(V^TV)}{\lambda_{min}(V^TV)}}\right\},
\]
where $\lambda_{\min}$ denotes the minimum eigenvalue of ${\cal M}:={\cal L}_d + {\cal P}$.
\end{theorem}}
\begin{proof}
Here, for convenience in robustness analysis, we assume that the
coupling strength $c=1$ in the controller (\ref{e4}). Under the conditions in Lemma \ref{l1}, the eigenvalues of ${\cal L}_d + {\cal P}$ are different from each other.
There exists a nonsingular matrix $V \in \mathbb{R}^{N\times N}$, $VV^{-1}=I_N$, such that
\begin{align}\label{e9}
{\cal L}_d + {\cal P}=V\Lambda V^{-1},
\end{align}
where $\Lambda = \mbox{diag}\{\lambda_1, \lambda_2, \ldots ,\lambda_N\}$, and $\lambda_i$ is the $i$-th real
eigenvalue of ${\cal L}_d + {\cal P}$.

Inserting (\ref{e9}) into (\ref{e7}), we can get
\[
G(s)=V\left[
        \begin{array}{cccc}
          G_1(s) &   &   &   \\
           &   G_2(s)   &   \\
            &   & \ddots &   \\
            &   &   & G_N(s) \\
        \end{array}
      \right]V^{-1}
\]
where
\[
G_i(s)=\frac{1}{\tau s^3+(1+\lambda_ik_a)s^2+\lambda_ik_vs+\lambda_ik_p},\quad i=1,\ldots,N.
\]

In the light of
\[
\gamma=\|G(s)\|_{{\cal H}_\infty} =\sup_{\omega}[\sigma_{\max}(G(j\omega))],
\]
and letting $\tilde{G}=\mbox{diag}\{G_i(s)\}_1^N$,
%\[
%\tilde{G}=\left[
%        \begin{array}{cccc}
%          G_1(s) &   &   &   \\
%           &   G_2(s)   &   \\
%            &   & \ddots &   \\
%            &   &   & G_N(s) \\
%        \end{array}
%      \right]
%\]
we have that
\begin{align*}
\gamma&=\|G(s)\|_{{\cal H}_\infty} =\sup_{\omega} \sqrt{\lambda_{max} [G^*(j\omega)G(j\omega)]}\\
&=\sup_{\omega} \sqrt{\lambda_{max} [V^{-T}\tilde{G}^*(j\omega)V^TV\tilde{G}(j\omega)V^{-1}]}\\
%&\leq \sqrt{\lambda_{max}(V^TV)} \sup_{\omega}  \sqrt{\lambda_{max}[V^{-T}\tilde{G}^*(j\omega)\tilde{G}(j\omega)V^{-1}]}\\
&\leq \sqrt{\lambda_{max}(V^TV)} \sup_{\omega} \max_i \sqrt{V^{-T}G_i^*(j\omega)G_i(j\omega)V^{-1}}\\%
&\leq \max_i \|G_i(s)\|_{{\cal H}_\infty}\sqrt{\frac{\lambda_{max}(V^TV)}{\lambda_{min}(V^TV)}}.
\end{align*}
where $^*(\cdot)$ means the complex conjugate operation. Furthermore, according
to $G_i(s)$, we can obtain
\[
\gamma \leq \frac{1}{\lambda_{\min}k_p}\sqrt{\frac{\lambda_{max}(V^TV)}{\lambda_{min}(V^TV)}}.
\]
{For the detail, see the proof of Theorem 1 in \cite{zllr2018}.}
And then,
\begin{align}\label{e10}
\gamma \leq \gamma_u:=\inf_{g_i,d_i,d_{ij}}\left\{\frac{1}{\lambda_{\min}k_p}\sqrt{\frac{\lambda_{max}(V^TV)}{\lambda_{min}(V^TV)}}\right\}.
\end{align}

$\hfill\blacksquare$
\end{proof}
From Remark \ref{r000}, we know that the matrix ${\cal M}$ can approximate to a diagonal matrix, if the parameters $d_i$ and $d_{i,j}$ are taken properly.
The matrix $V$ approximates to $\bar{\kappa} I$ ($\bar{\kappa}$ is a constant).
Further, we have

\begin{align*}%\label{e10}
\gamma \leq \gamma_u:=\inf_{g_i,d_i,d_{ij}}\left\{\frac{1}{\lambda_{\min}k_p}\bar{\kappa}\right\}.
\end{align*}
{\begin{remark}
Theorem \ref{t1} shows an implicit relation between the proposed
robustness measure, i.e. $\gamma$-gain, and the introduced weighting parameters. In the proof, the key point lies in the diagonalization of ${\cal M}$ matrix,
thus it results in the decoupling of the transfer function $G(s)$.
{Therefore, it only needs to analyze the ${H}_\infty$
norm of a single transfer function, and this transfer function can be modified by adjusting the eigenvalue in the light of
the parameters $g_i, d_i$ and $d_{ij}$ in the controller (\ref{e4}).} $\hfill\blacksquare$
\end{remark}}

{\begin{remark}\label{r2} To achieve the better robustness performance for a mass of platoons,
(\ref{e10}) implies that we must select appropriate controller parameters $g_i, d_i$ and $d_{ij}$  associated with a larger
$\lambda_{\min}$.
A practical choice is to make the center of the leftmost Ger\v{s}gorin disc be far away from the origin as possible,
and then its radius as small as possible (See Fig. \ref{f22}). It means that  the minimum eigenvalue $\lambda_{\min}$ is enlarged.
$\hfill\blacksquare$
\end{remark}}

Different from the results in \cite{zllr2018} (see Corollaries 1 and 2), Remark \ref{r2} shows that under the presented controller (\ref{e4}),
the robustness performance for a mass of platoons does depend on the controller parameters $g_i, d_i$ and $d_{ij}$,
and is no longer limited to the total number $N$ of vehicles and the number of followers that are connected to the
leader (structure).
It will be a more flexible way to improve the robustness performance for a
mass of platoons without changing the structure of the platoon and the number of vehicles.
% needed in second column of first page if using \IEEEpubid
%\IEEEpubidadjcol

\section{Distributed ${ H}_\infty$ Controller Synthesis
Of the Platooning}\label{s4}
Here, we introduce a modified distributed controller for a platooning with the guaranteed ${ H}_\infty$ performance.
To satisfy $\gamma= \|G(s)\|_{{ H}_\infty} < \gamma_d$ with a given desired $\gamma_d$-gain, the feedback gains ($(k = [k_p\quad k_v\quad k_a]^T)$) and
coupling strength $c$ need to be designed.

Based on the decoupling technique in the proof of Theorem \ref{t1}, the distributed ${H}_\infty$ control problem is converted into a set of ${ H}_\infty$ control
of independent systems sharing the same dimension with a
single vehicle (See Theorem 3 in \cite{ldc2011}). %Lemma \ref{l2}).

%Lemma \ref{l2}
The collective behavior of a platoon system (\ref{e6}) is decoupled
into a set of $N$ individual subsystems (\ref{e12})
\begin{align}\label{e12}
\dot{\bar{x}}_i &= (A - c\lambda_i B_1k^T )  \bar{x}_i + B_2 \bar{w}_i ,\nonumber\\
\bar{y}_i &= C_1 \bar{x}_i ,\quad i = 1,\ldots, N,
\end{align}
in the light of Theorem 3 in \cite{ldc2011}.
%The influence of the introduced weighting parameters shows the fact that the decoupled systems are modified by the eigenvalues of ${\cal L}_d+{\cal P}$. Also,
The distributed ${ H}_\infty$ control
problem of a platooning (\ref{e6}) is equivalently reduced to a set of ${ H}_\infty$ control problems of subsystems (a single vehicle) with the
same dimension $3 \times 3$.
The computational complexity is dramatically reduced due to the less dimension of systems, for the controller synthesis.

Next we synthesize the distributed ${\cal H}_\infty$ controller by the following theorem.
\begin{theorem}\label{t2}\cite{zllr2018}
Consider a homogeneous platooning (\ref{e6}) with the directed
topology. For any desired $\gamma_d > 0$, we have
%the following statement:
%\begin{enumerate}
    %For the closed-loop platoon system,
$\|G(s)\|_{{\cal H}_\infty} < \gamma_d$, if the feedback gains are chosen
as $k^T = \frac{1}{2}B^T_1 Q^{-1}$, and the coupling strength satisfies
\begin{align} \label{e13}
c\geq \frac{\sqrt{\alpha}}{\lambda_{\min}},
\end{align}
where $Q >0$ and $\alpha > 0$ are the feasible solutions to the linear matrix inequality (\ref{e11}) % in  Lemma 5 . Lemma \ref{l3}.
 \begin{align}\label{e11}
\left[
  \begin{array}{ccc}
    AQ+QA^T-\alpha B_1B_1^T & B_2 & QC_1^T \\
    B_2^T & -\gamma_d^2 & 0 \\
    C_1Q & 0 & -1 \\
  \end{array}
\right]<0.
\end{align}
% \item There always exists a feasible solution: $Q>0$ and $\alpha > 0$ to (\ref{e11})  %Lemma \ref{l3};
%\end{enumerate}
\end{theorem}

To design a distributed ${\cal H}_\infty$ controller, the feedback gains $k^T$ can be obtained by solving LMI (\ref{e11}), and
the coupling strength $c$ can be adjusted to satisfy the condition (\ref{e13}).
In addition, it retains the advantage of the computational complexity in \cite{zllr2018}.
\begin{remark}\label{r3}
It is worth mentioning that Theorem \ref{t2} implies that for a platoon system,
a distributed controller (\ref{e4}) can be synthesized to satisfy any
given ${\cal H}_\infty$ performance.
The new controller can circumvent the high-gain case mentioned in \cite{zllr2018},
if only the matrix ${\cal L}_d + {\cal P}$ associated with the directed graph has a
larger $\lambda_{\min}$.
In fact, we can do it by moving the leftmost Ger\v{s}gorin disc away from the origin as far as possible,
in the light of adjusting the controller parameters $g_i, d_i$ and $d_{ij}$.
{It not only improves the robustness performance for a given controller (See the expression (\ref{e10})),
but also reduces the value of  coupling strength $c$ for a given ${ H}_\infty$ performance (See the expression (\ref{e13})).}
$\hfill\blacksquare$
\end{remark}

\section{Numerical Simulation}\label{s5}
In this section, one numerical experiment with passenger cars are used to verify the effectiveness of the proposed method for directed graphs. %, compared with \cite{zllr2018}.
%First, the scaling trend of $\gamma$-gain under the given directed graph is demonstrated in Table \ref{tab:cspi}.
Next the synthesis of a distributed ${H}_\infty$ controller is also illustrated. %In addition, simulations with a realistic model
Here we consider that the generalized directed topology with $8$ vehicles and one leader vehicle $0$ which is shown in Fig. \ref{f1}.

\begin{figure}[H]
\begin{center}
  % Requires \usepackage{graphicx}
  \includegraphics[width=3.5in]{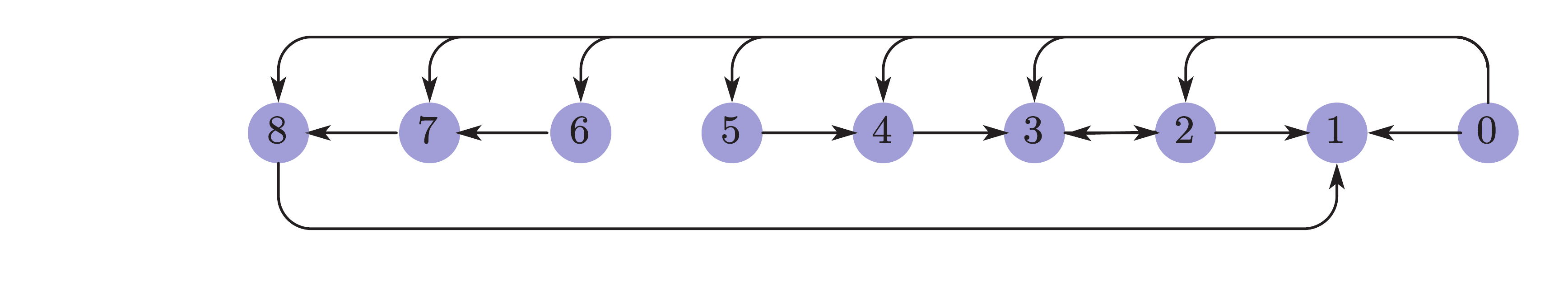}
  \caption{Directed communication topology with $8$ vehicles and one leader vehicle}\label{f1}
  \end{center}
\end{figure}
{\subsection{Validation of superiority}}
We can design a distributed ${ H}_\infty$ controller in the light of Theorem \ref{t2}, using the same parameters as that in \cite{zllr2018}
(i.e., $\tau=0.5 s$, $\alpha=1.968$ and $k=[2.122\quad 3.425\quad 2.501]^T$ under the desired performance $\gamma_d=1$). We also implement these controllers for the platooning, like the scenario in \cite{zllr2018}: the zero initial tracking state
errors for the platooning, a constant reference speed for the leader ($v_0 = 20 m/s$), and the same external
disturbances for each node:
\begin{align*}
w_i(t)=\left\{
\begin{aligned}
0\qquad \qquad \qquad \qquad &0<t<5s \\
Q\sin(\frac{2\pi}{5}(t-5))\qquad &5s\leq t< 10s \\
0 \qquad \qquad \qquad \qquad &10s\leq t
\end{aligned}
\right.
\end{align*}
where $Q=10$. {The disturbance $w_i$ is a generalized expression for all possible disturbance (such as aerodynamic drag, frictional drag and slope, etc.), hence its physical significance hinges on the application scenarios.}

Under the controller (\ref{e4}), the matrix ${\cal M}$ corresponding to the topology in Fig. \ref{f1} is
\[
{\cal M}=\left[
  \begin{array}{cccccccc}
    \bar{g}_1 & -1 & 0 & 0 & 0 & 0 & 0 & -1 \\
    0 & \bar{g}_2  & -1 & 0 & 0 & 0 & 0 & 0 \\
    0 & -1 & \bar{g}_3 & -1 & 0 & -1 & 0 & 0 \\
    0 & 0 & 0 & \bar{g}_4 & -1 & 0 & 0 & 0 \\
    0 & 0 & 0 & 0 & g_5 & 0 & 0 & 0 \\
    0 & 0 & 0 & 0 & 0 & g_6 & 0 & 0 \\
    0 & 0 & 0 & 0 & 0 & -1 & \bar{g}_7 & 0 \\
    0 & 0 & 0 & 0 & 0 & 0 & -1 & \bar{g}_8 \\
  \end{array}
\right],
\]
with $d_{ij}=1\ (j\in \mathbb{N}_i), \bar{g}_1=2d_1+g_1, \bar{g}_3=3d_3+g_3, \bar{g}_i=d_i+g_i, i=\{2,4,7,8\}$.

First of all, we chose the parameters of Test (a) in Table \ref{tab:cp} for the ${ H}_\infty$ controller.
And then we can calculate numerically $\lambda_{\min}({\cal L}_d+{\cal P})$ and the coupling strength $c$,
the error amplification $\frac{\|Y (t)\|_{{\cal L}_2}}{\|W(t)\|_{{\cal L}_2}}$ (i.e., $\gamma$-gain $\gamma$) in time-domain,
%and its upper bound $\gamma_u$,
see Test (a) in Table \ref{tab:cspi}.
The curves of spacing errors are shown in Fig. \ref{fig1}.
It shows that the platooning is able to regain stability after receiving the disturbance from $5s$ to $10s$,
under the ${ H}_\infty$ controller which is synthesized by Theorem \ref{t2}.
{This simulation verifies the effectiveness of our controller for the directed graph.
While the work \cite{zllr2018} can not deal with the this case for the directed graph.}
\begin{figure}[H]
  \centering
%\hspace{-10mm}
  \includegraphics[width=3.5in]{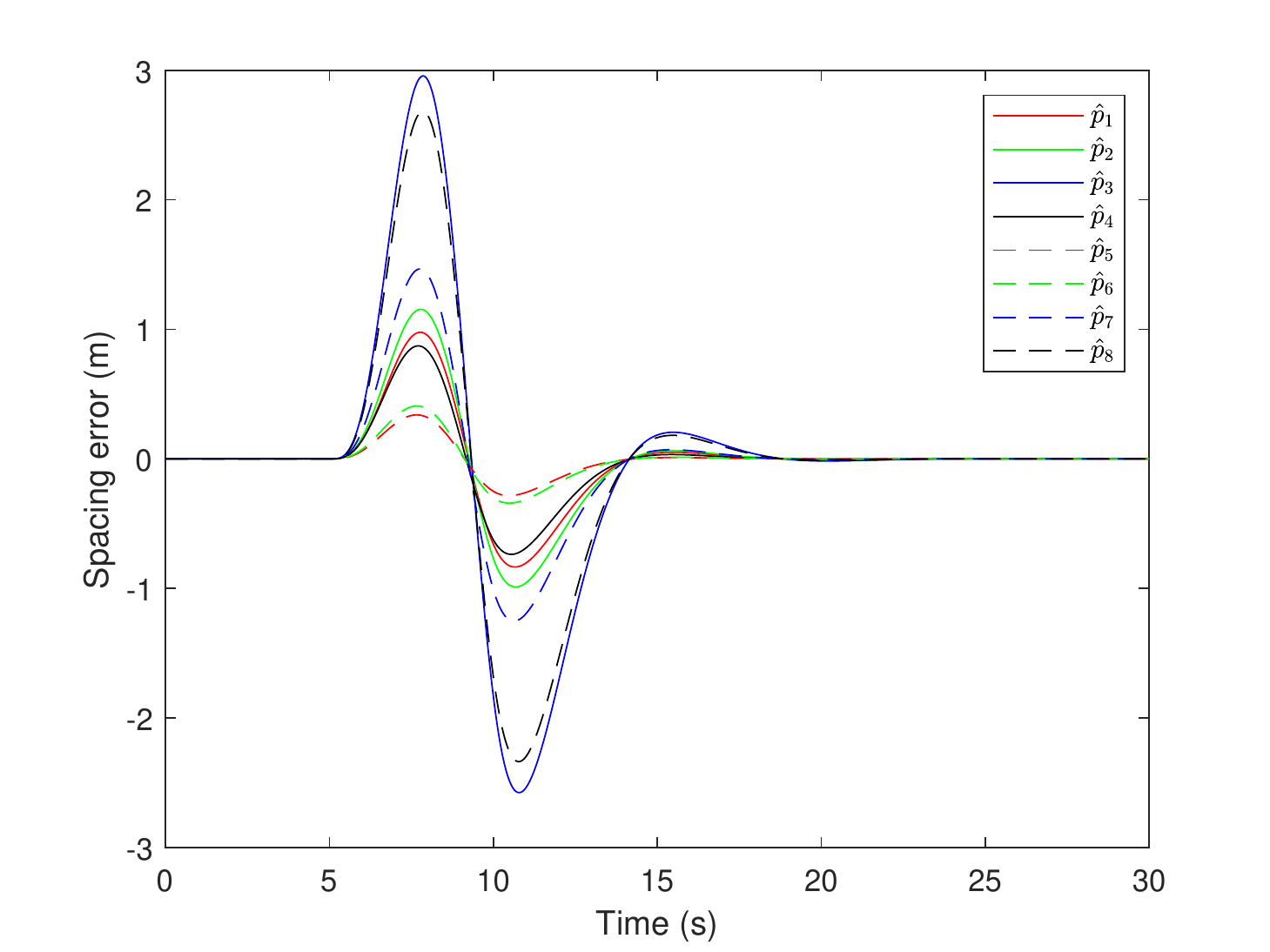}
  \caption{Spacing errors for platooning under distributed ${ H}_\infty$ controllers with $\lambda_{\min}({\cal L}_d+{\cal P})=2.1$.}\label{fig1}
\end{figure}

{Next we change the parameters as Test (b) shown in Table \ref{tab:cp}.
The corresponding results are obtained in the same way by the numerical calculation, as Test (b) shown in Table \ref{tab:cspi}.
Meanwhile, the corresponding profiles of spacing errors are shown in Fig. \ref{fig2}. Obviously, the worst spacing error became
smaller under the same external
disturbances. Thus it has the better robustness with the larger $\lambda_{\min}({\cal L}_d+{\cal P})$.
In addition, the coupling strength $c$ does not have the high gain, and the robustness is also improved, compared with that in \cite{zllr2018}.}
Note that the sign `--' means that the data is not provided in \cite{zllr2018}.}
 \begin{figure}[H]
  \centering
%\hspace{-10mm}
  \includegraphics[width=3.5in]{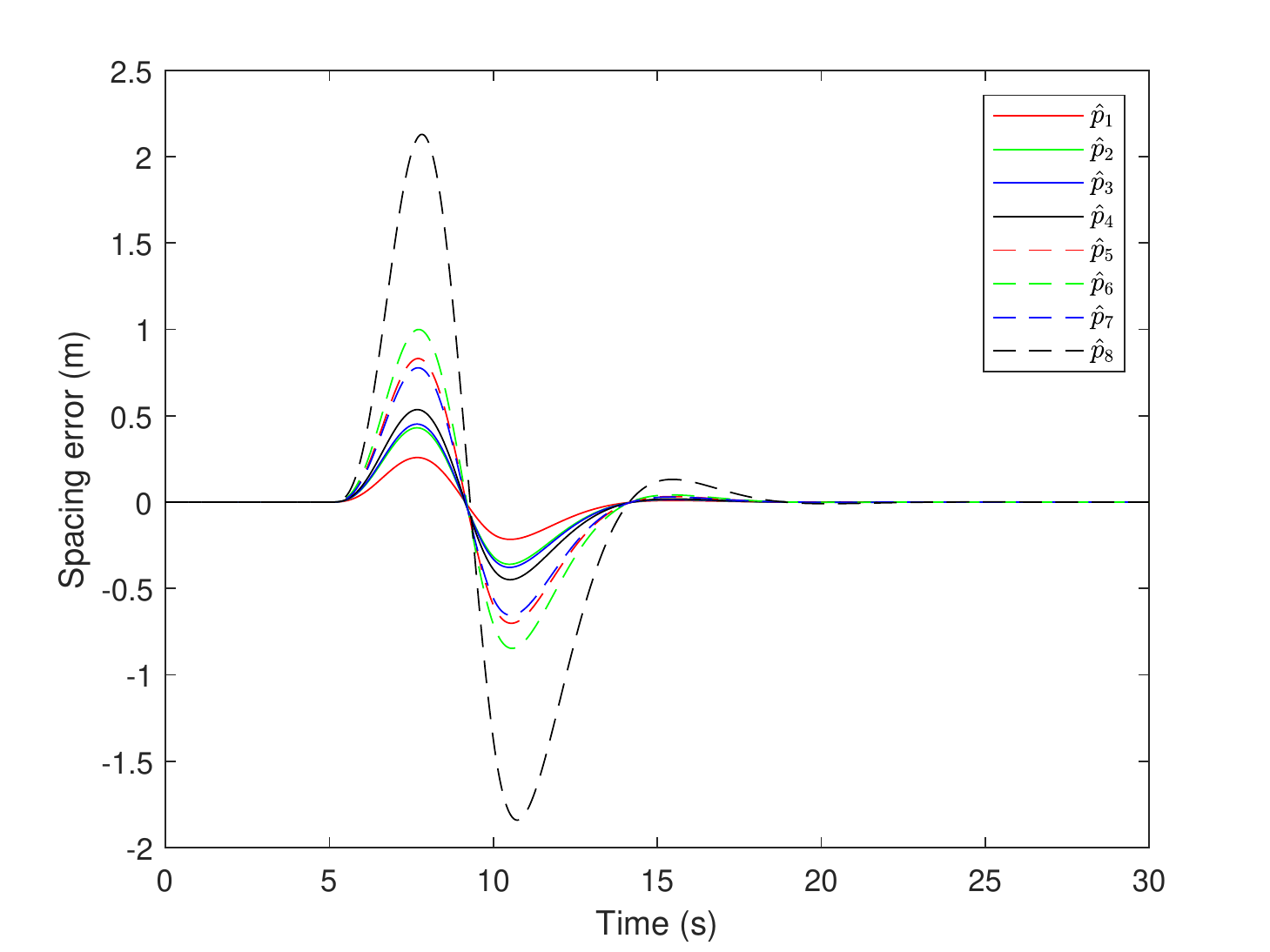}
  \caption{Spacing errors for platooning under distributed ${ H}_\infty$ controllers with $\lambda_{\min}({\cal L}_d+{\cal P})=5.1$.}\label{fig2}
\end{figure}
\vspace{-.5cm}
\begin {table}[H]
\caption {Controller parameters} \label{tab:cp}
\begin{center}
\begin{tabular}{cccccccccc}
\toprule
Test &$d_1$ & $d_2$ & $d_3$ & $d_4$ & $d_7$& $d_8$ & $g_5,g_6$ & $g_i$ ($i\neq 5,6$) &$d_{ij}$ \\
\midrule
%(a)&60& 140 & 150 & 160 & 110 & 130 & 155 & $\kappa$ & 1 \\
%\midrule
%(a)&8& 9 & 2 & 10 & 7 & 4 & 5 & 4 & 1\\
(a)&4& 6 & 1 & 5 & 3 & 2 & 12,  10 &0.1 &1\\
\midrule
(b)&24& 24 & 12 &20& 7 & 14 & 12,  10 & 0.1&1\\
\bottomrule
\end{tabular}
\end{center}
\end{table}
\vspace{-.5cm}
\begin {table}[H]
\caption {Calculating coupling strength and some performance indices} \label{tab:cspi}
\begin{center}
\begin{tabular}{cccccccccc}
\toprule
Test & $\lambda_{\min}({\cal L}_d+{\cal P})$ & $c$ & $\gamma$\\%& $\gamma_u$
\midrule
(a)&2.1 & 0.6680 & 0.4501 \\%& 2.7647\\
\midrule
(b)&5.1 & 0.2751 &  0.2996\\% & 0.5598\\
\midrule
\cite{zllr2018}&0.0557&35.33&--\\%&--\\
\bottomrule
\end{tabular}
\end{center}
\end{table}

\subsection{Influence of disturbance for the platoon system}
The spacing errors of the platoon are also mainly dependent on the amplitude (strength) of disturbance $w_i$ under the given ${ H}_\infty$ controllers. It is well-known that the errors will generally become larger as the amplitude of $w_i$ grows.
{In the following, we present the simulation examples to verify this point in the light of the controller which is used in Fig \ref{fig1}.
We take $Q=30, 5s\leq t\leq10s$, then the spacing errors lie in the interval $(0.4, 6.4)$ (meters), which are shown in Fig. \ref{fig5}.
 \begin{figure}[H]
  \centering
%\hspace{-10mm}
  \includegraphics[width=3.5in]{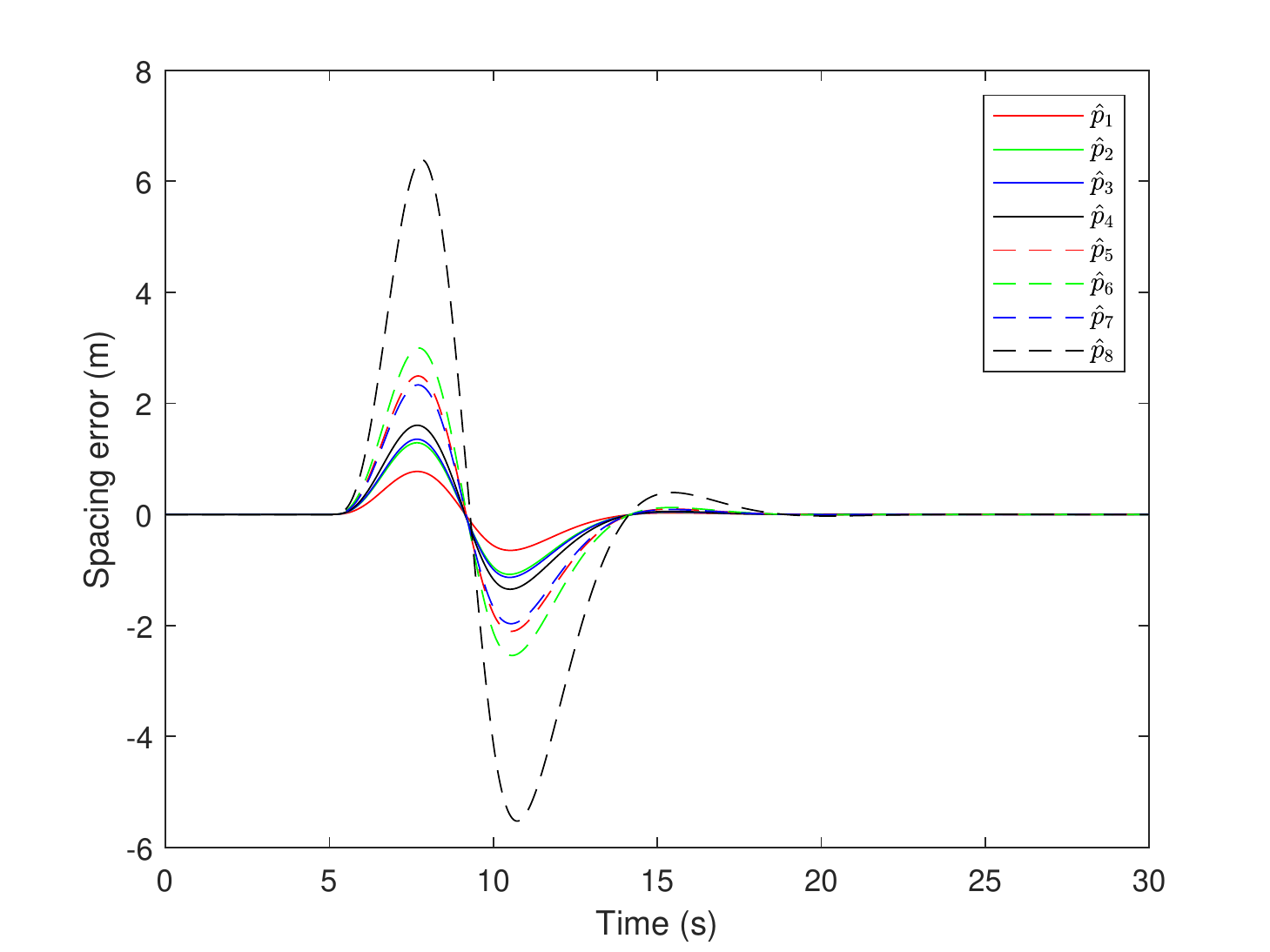}
{ \caption{Spacing errors for platooning under the same controllers in Fig \ref{fig1} with $Q=30$.}\label{fig5}}
\end{figure}}
{Obviously, the errors are larger than ones in Fig.\ref{fig1} where the errors are in the interval $(0.3,2.9)$ (meters). Thus we verify that the errors will become large with the increasing of strength of disturbance.}

{\subsection{A simulation besed on NGSIM dataset}
Here, the Next Generation Simulation (NGSIM) dataset is used to verify the effectiveness of our method.
In the dataset, the real-world trajectory information (including lane-specific location, velocity and
acceleration, collected every one-tenth of a second) of the entire traffic flow were available.
We screen out the information (time, longitudinal location, velocity and acceleration,
referring to \cite{Sun2021}) of vehicle (No.4) in the I-80 dataset as the state of the leading car (the node $0$ in Fig. \ref{f1}).}

{Due to the NGSIM trajectory data exhibit errors and noises that need to be filtered out,
we adopt the locally weighted regression method \cite{Toledo1999} to get their smooth and continuous curve, see Fig. \ref{fig6}.
\begin{figure}[H]
  \centering
%\hspace{-10mm}
  \includegraphics[width=3.3in]{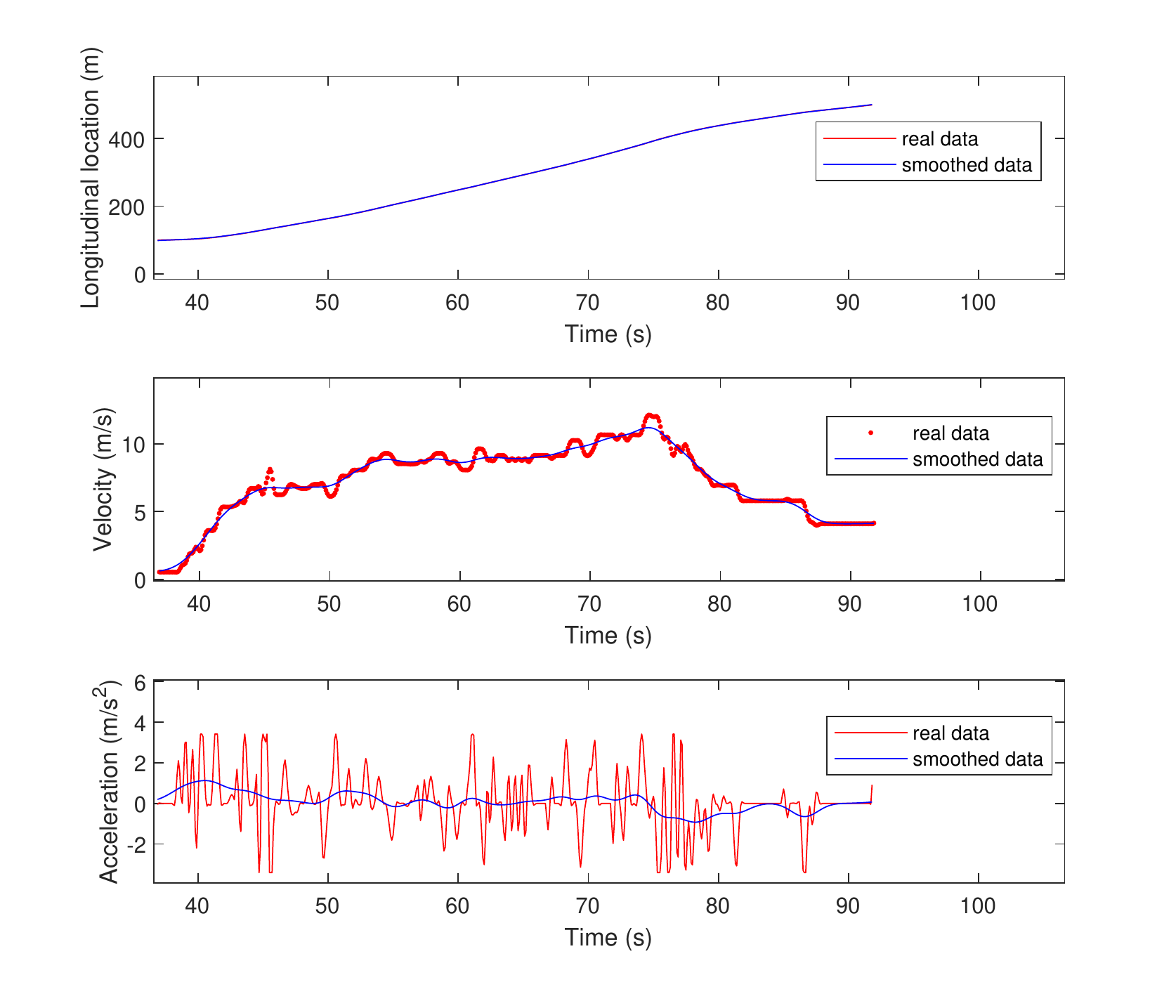}
 \caption{The smoothed curve of the raw data.}\label{fig6}
\end{figure}}

{Next, we use the designed controller in Test (a) of Table \ref{tab:cp} to steer the vehicles to
reach a consensus on speed and keep the desired spacing ($|d_{i,i-1}|=4.5m$) over the communication topology in Fig. \ref{f1}.
To verify the robustness of the vehicle platoon under the controller, we exert a disturbance (aerodynamic drag $c_2v_i^2$, where $c_2=0.5$) into the platoon system in interval $(57s, 63s)$.}

\begin{figure}[H]
  \centering
%\hspace{-10mm}
  \includegraphics[width=3.5in]{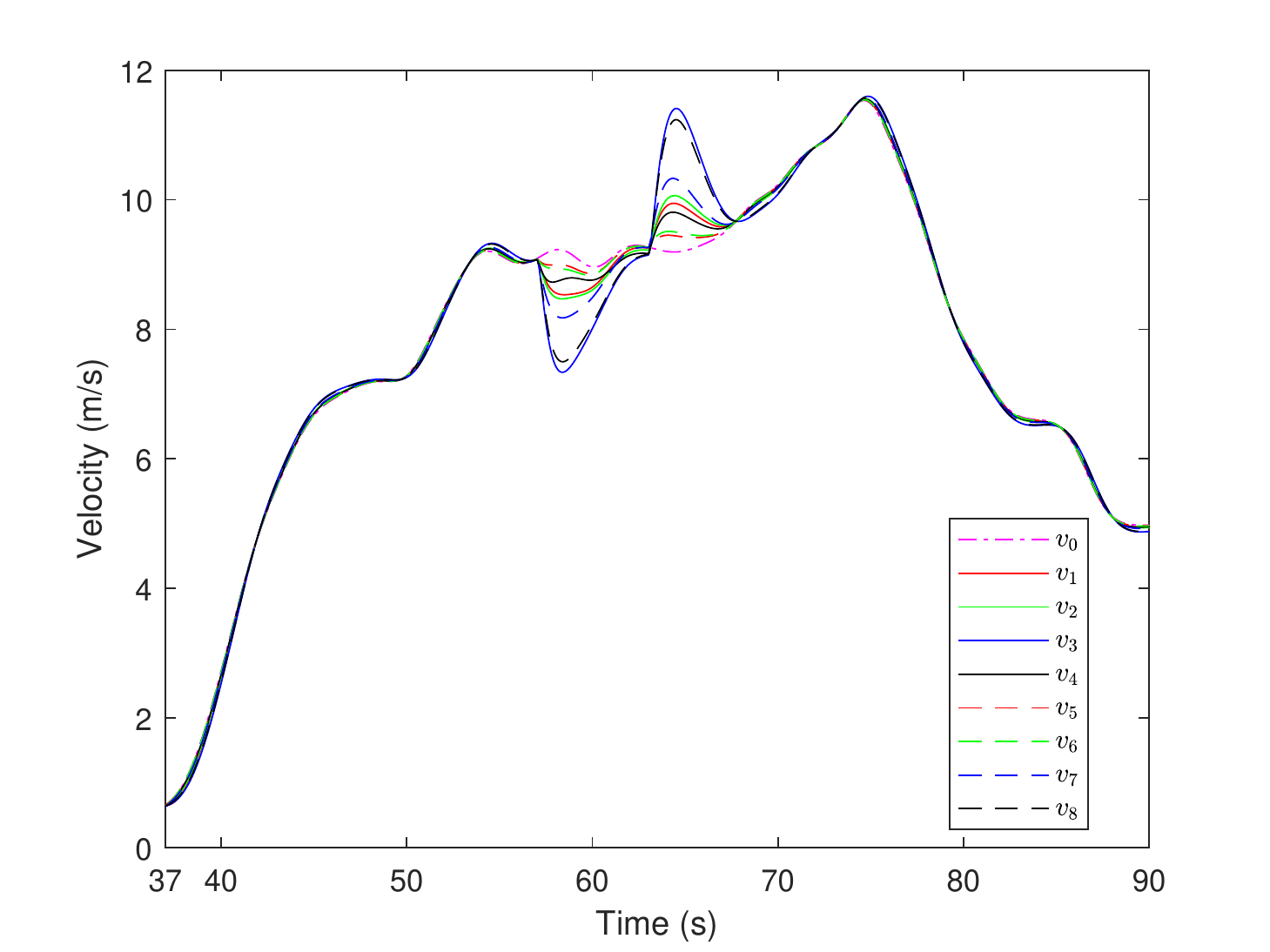}
 \caption{The trajectories of velocity tracking.}\label{f7}
\end{figure}

{From Fig. \ref{f7}, except for the time interval $(57s, 71s)$, the effect of tracking is very good.
From 57s to 63s, the velocity of the following cars slow down due to the effect of aerodynamic drag. And then, after 71s, it reaches a consensus on $v_0$ again.
The phase from 63s to 71s is the transient process.}

\begin{figure}[H]
  \centering
%\hspace{-10mm}
  \includegraphics[width=3.5in]{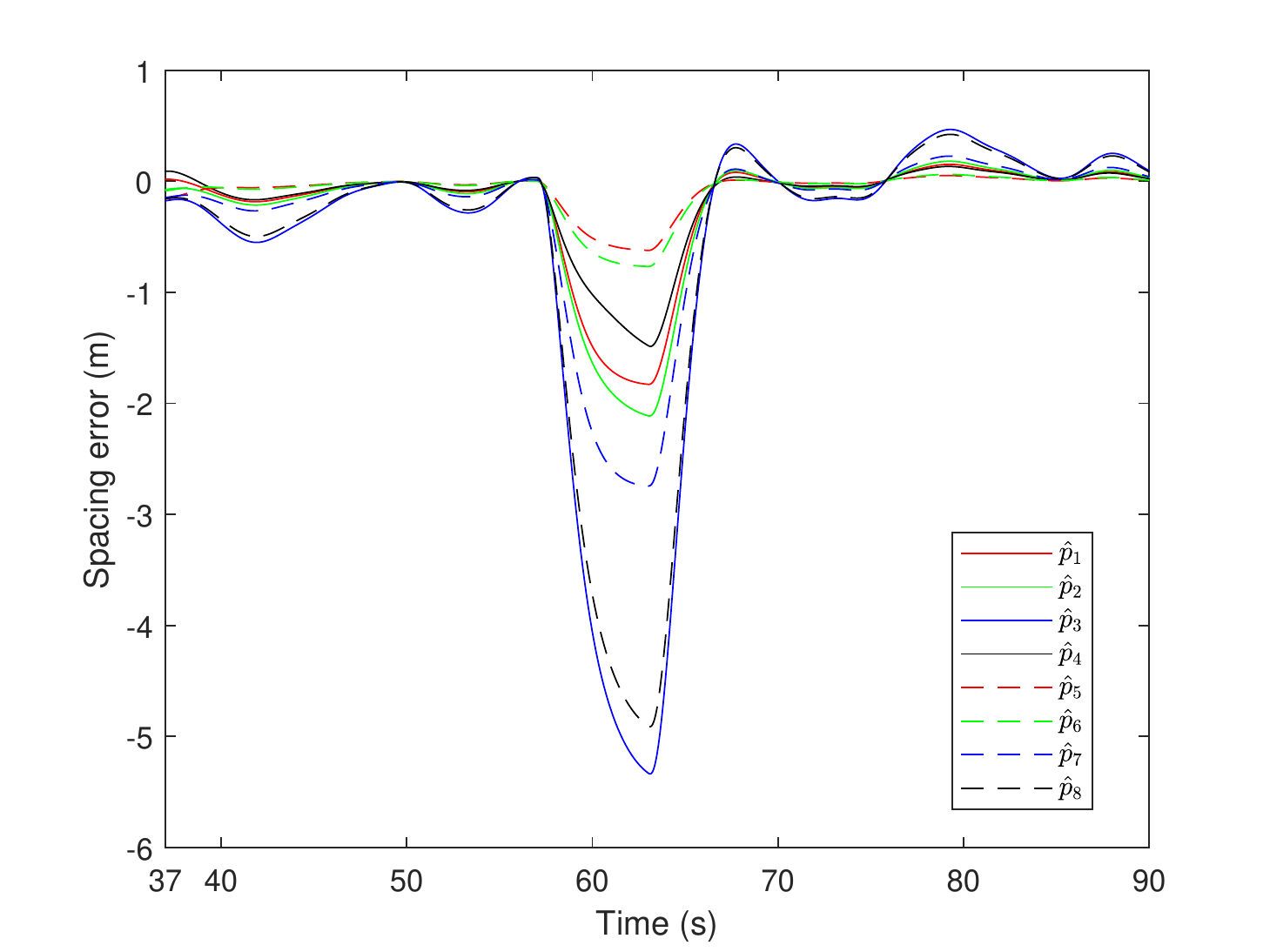}
 \caption{The trajectories of spacing error between the front vehicle and rear vehicle.}\label{f9}
\end{figure}

{From Fig. \ref{f9}, the spacing error of every following vehicle has changed when the disturbance is exerted.
The vehicles $3$ and $8$ are the most obvious ($|\hat{p}_3|=5.3m$ and $|\hat{p}_8|=4.9m$). When the disturbance vanishes, the spacing errors go to zero again, that is, the platoon system is back to the desired state.}

\begin{figure}[H]
  \centering
  \includegraphics[width=3.5in]{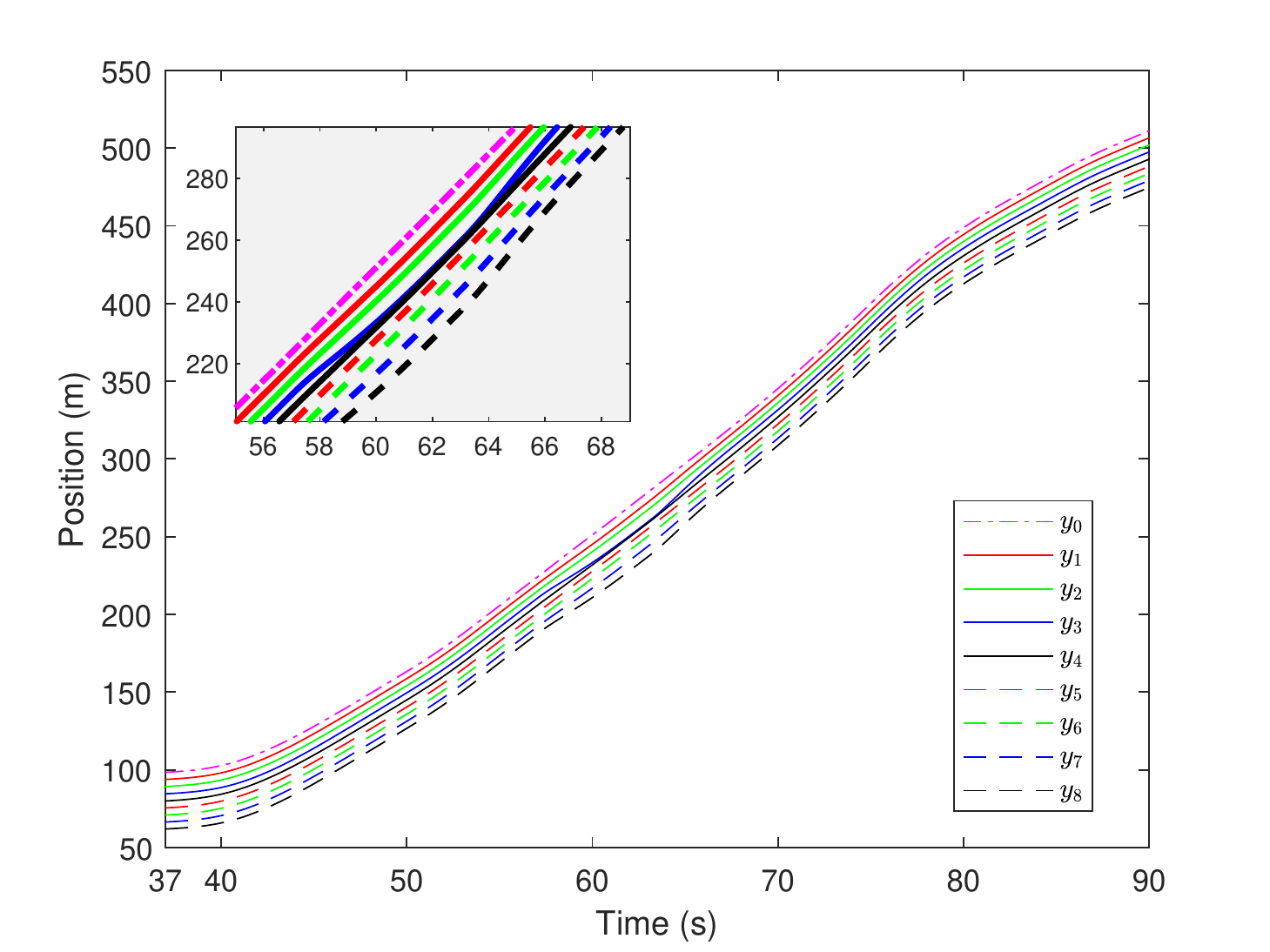}
 \caption{The trajectories of longitudinal location.}\label{f8}
\end{figure}

Above result is also verified by the trajectories of longitudinal location in Fig. \ref{f8}, the speed of vehicles $3$ and $8$ drop off notably when the aerodynamic drag occurs.
And then the vehicle $3$ nearly has a collision with the vehicle $4$.
It is not what we would like to see, since there is a high probability of traffic accidents. Note that although the spacing error of the vehicle $3$ is more than the
desired spacing ($4.5m$), they do not collide, since the vehicle $4$ also
slows down by the effect of the aerodynamic drag.

\begin{figure}[H]
  \centering
%\hspace{-10mm}
  \includegraphics[width=3.5in]{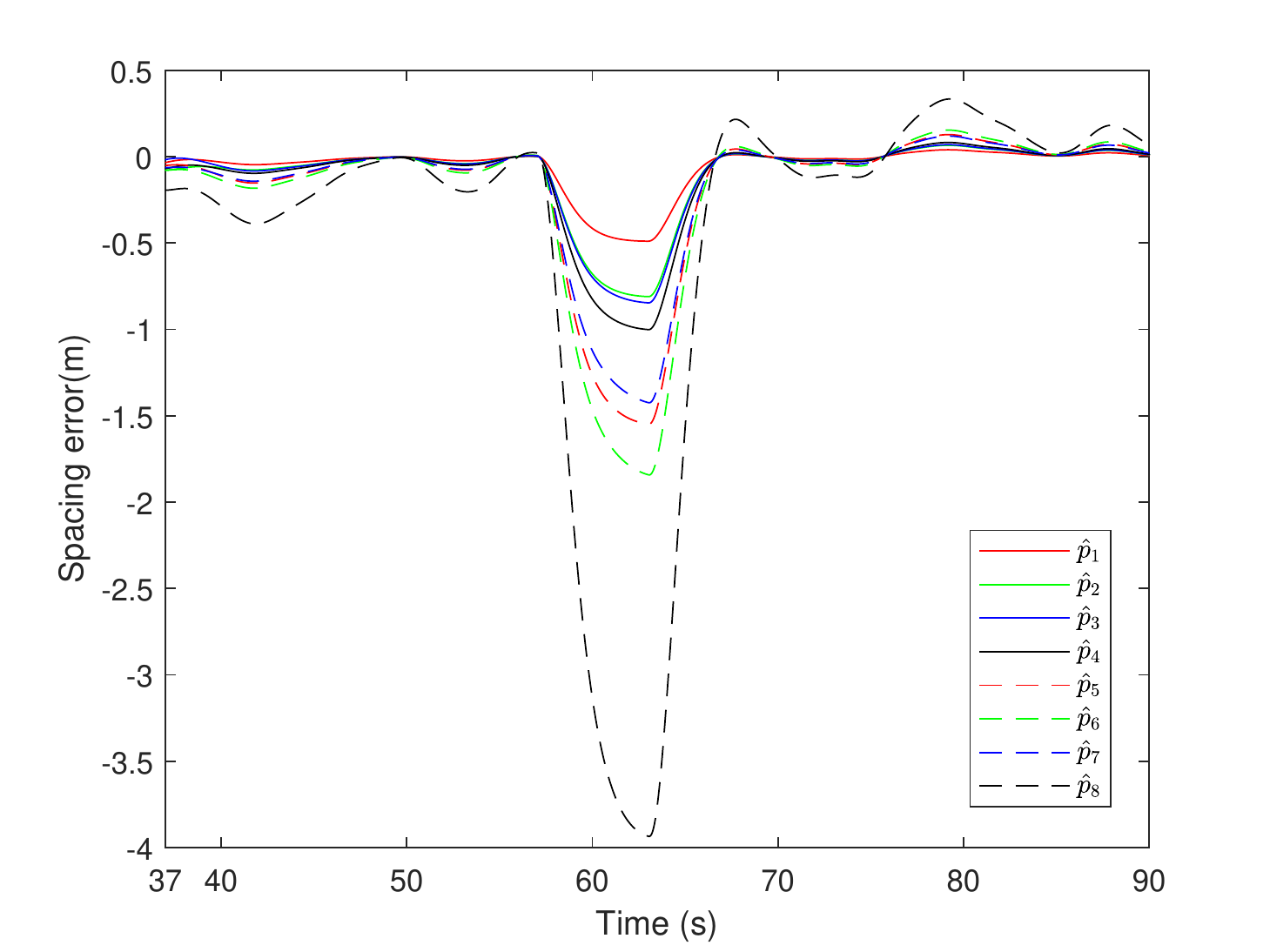}
 \caption{The trajectories of spacing error between the front vehicle and rear vehicle with Test (b).}\label{f10}
\end{figure}

\begin{figure}[H]
  \centering
%\hspace{-10mm}
  \includegraphics[width=3.5in]{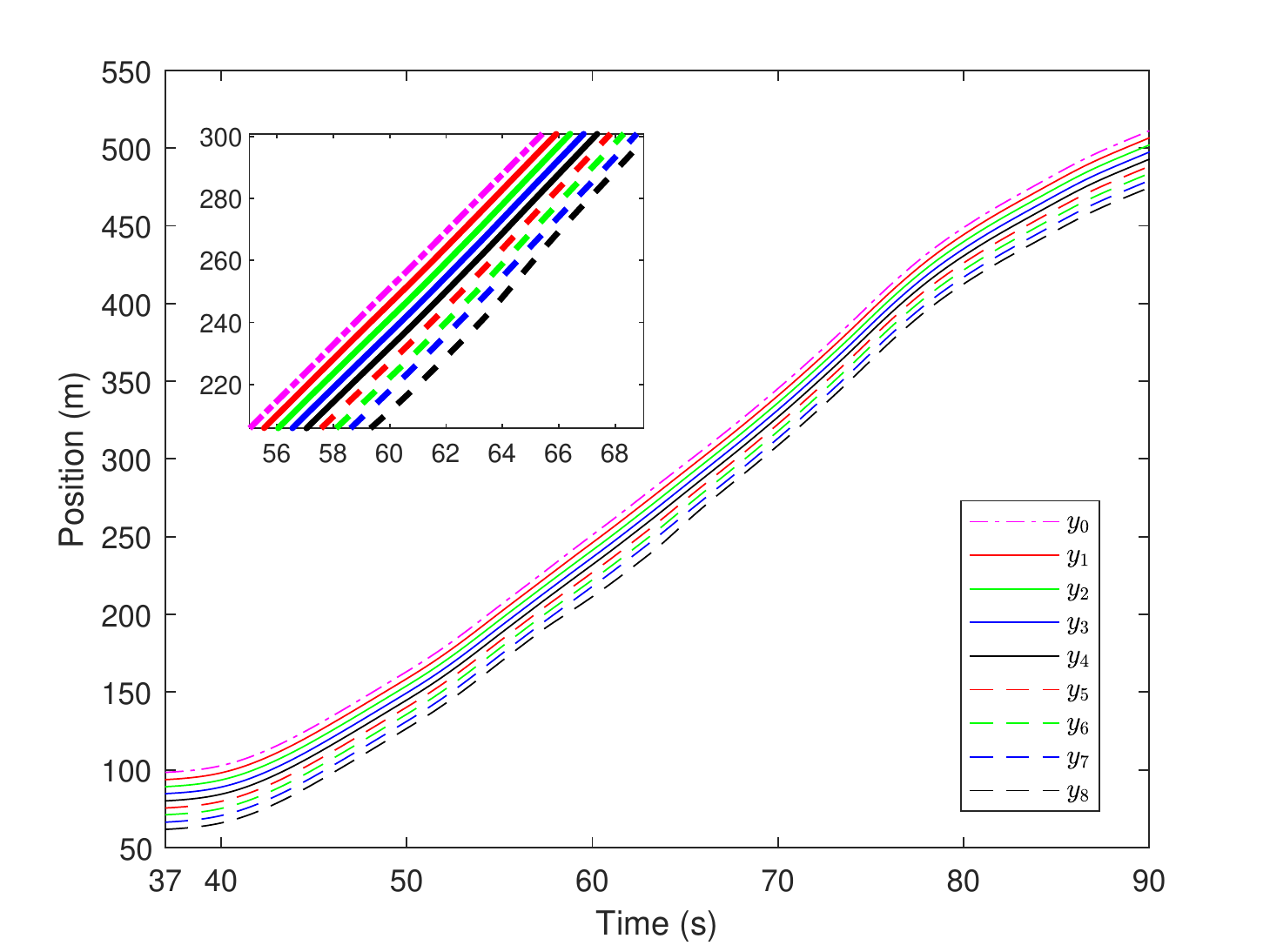}
 \caption{The trajectories of longitudinal location with Test (b).}\label{f11}
\end{figure}

To avoid the collision between the front vehicle and the rear one,
the robustness of the platoon system is improved by our method,
in the light of the controller parameters in Test (b) of Table \ref{tab:cp}.
Apparently, the worst of spacing error is improved, from Fig. \ref{f10}, that is, the spacing error is reduced to $3.9m$ from $5.3m$ compared with that in Fig. \ref{f9}.
Moreover, Fig. \ref{f11} shows that no vehicle is at risk of collision when the disturbance occurs.
The above simulation verifies the effectiveness of our method for improving the robustness of platoon system by NGSIM data.

\section{Conclusion}\label{s6}
In this paper, a new control structure is presented to study the robustness and distributed
${ H}_\infty$ controller synthesis for a platooning of connected vehicles with the directed topology.
Compared with the problem of the conventional undirected topology \cite{zllr2018},
we solve the problem on the directed topology, which is more challenging.
This paper features itself in two aspects:
first, the robustness performance of the platoon control can be improved by choosing the suitable controller parameters, and
it no longer depends on the number of vehicles $N$ and information topology;
second, we circumvent the high-gain problem for the synthesis of distributed ${ H}_\infty$ controller.
Namely the coupling strength $c$ is not very large such that the new control structure is practical in implementation.
In the future, we will try to study the string stability problem by the proposed controller. We think it is feasible, since our method has generality for platoon systems.

\ifCLASSOPTIONcaptionsoff
  \newpage
\fi

% trigger a \newpage just before the given reference
% number - used to balance the columns on the last page
% adjust value as needed - may need to be readjusted if
% the document is modified later
%\IEEEtriggeratref{8}
% The "triggered" command can be changed if desired:
%\IEEEtriggercmd{\enlargethispage{-5in}}

% references section

% can use a bibliography generated by BibTeX as a .bbl file
% BibTeX documentation can be easily obtained at:
% http://mirror.ctan.org/biblio/bibtex/contrib/doc/
% The IEEEtran BibTeX style support page is at:
% http://www.michaelshell.org/tex/ieeetran/bibtex/
%\bibliographystyle{IEEEtran}
% argument is your BibTeX string definitions and bibliography database(s)
%\bibliography{IEEEabrv,../bib/paper}
%
% <OR> manually copy in the resultant .bbl file
% set second argument of \begin to the number of references
% (used to reserve space for the reference number labels box)

\bibliographystyle{IEEEtran}
\bibliography{mydatabase}

%\begin{thebibliography}{1}

%\bibitem{IEEEhowto:kopka}
%H.~Kopka and P.~W. Daly, \emph{A Guide to \LaTeX}, 3rd~ed.\hskip 1em plus
%  0.5em minus 0.4em\relax Harlow, England: Addison-Wesley, 1999.

%\end{thebibliography}

% biography section
%
% If you have an EPS/PDF photo (graphicx package needed) extra braces are
% needed around the contents of the optional argument to biography to prevent
% the LaTeX parser from getting confused when it sees the complicated
% \includegraphics command within an optional argument. (You could create
% your own custom macro containing the \includegraphics command to make things
% simpler here.)
\vspace{-.5in}
\end{document}